\documentstyle[amscd,amssymb,amsthm,verbatim,epsf]{amsart}

\begin{document}
\theoremstyle{plain}
\newtheorem{Thm}{Theorem}
\newtheorem{Cor}{Corollary}
\newtheorem{Con}{Conjecture}
\newtheorem{Main}{Main Theorem}
\newtheorem{Lem}{Lemma}
\newtheorem{Prop}{Proposition}
\newtheorem{Exp}{Example}
\theoremstyle{definition}
\newtheorem{Def}{Definition}
\newtheorem{Note}{Note}

\theoremstyle{remark}
\newtheorem{notation}{Notation}
\renewcommand{\thenotation}{}

\errorcontextlines=0
\numberwithin{equation}{section}

\renewcommand{\rm}{\normalshape}%
\title [Abel Continuity]%
   {Abel Continuity}
\author{H\"usey\.{I}n \c{C}akall\i*, and Mehmet Albayrak** \\ *Maltepe University, Faculty of Science and letters, Department of Mathematics, Istanbul, Turkey Phone:(+90216)6261050, fax:(+90216)6261113\\{**}Sakarya University, Adapazari, Turkey Phone:(+90264)2955946, fax:(+90264)2955950}

\address{H\"usey\.{I}n \c{C}akall\i \\
          Maltepe University, Department of Mathematics,  \; \;  \; \; \; \;  \; \; \; \; \; \; Faculty of science and letters, Marmara E\u{g}\.{I}t\.{I}m K\"oy\"u, TR 34857, Maltepe, \.{I}stanbul-Turkey \; \; \; \; \; Phone:(+90216)6261050 ext:1206, \;  fax:(+90216)6261113}

\email{hcakalli@@maltepe.edu.tr; hcakalli@@gmail.com}

\address{Mehmet ALBAYRAK \\
           Sakarya University, Institute of Natural Sciences, Department of Mathematics, Adapazari, Turkey Phone:(+90264)2955946, fax:(+90264)2955950}

\email{mehmetalbayrak12@@gmail.com}


\date{\today}
\begin{abstract}
A sequence $\textbf{p}=(p_{n})$ of real numbers is called Abel
convergent to $\ell$ if
the series $\Sigma_{k=0}^{\infty}p_{k}x^{k}$ is convergent for $0\leq x<1$ and
\[\lim_{x \to 1^{-}}(1-x) \sum_{k=0}^{\infty}p_{k}x^{k}=\ell.\]
We introduce a concept of Abel continuity in the sense that a function $f$ defined on a subset of $\Re$, the set of real numbers, is Abel continuous if it preserves  Abel convergent sequences, i.e. $(f(p_{n}))$ is an Abel convergent sequence whenever $(p_{n})$ is. A new type of compactness, namely Abel sequential compactness is also introduced, and interesting theorems related to this kind of compactnes, Abel continuity, statistical continuity, lacunary statistical continuity, ordinary continuity, and uniform continuity are obtained.

\end{abstract}

\maketitle

\section{Introduction}

The concept of continuity and any concept involving continuity play very important role not only in pure mathematics but also in other branches of sciences
involving mathematics especially in computer sciences, information theory, and dynamical systems.

It is well known that a real function $f$ is continuous if and only if $(f(p_n))$ is convergent
whenever $(p_n)$ is, and a subset $E$ of $\Re$ is compact if any sequence $(p_n)$ of points
in $E$ has a convergent subsequence whose limit is in $E$ where $\Re$ is the set of real
numbers. Combining this idea with Abel convergent sequences we introduce new concepts, namely, Abel continuity, and Abel sequential compactness.

The purpose of this paper is to study the concept of Abel continuity for real functions and show that this kind of continuity implies ordinary continuity, and present theorems related to this type of continuity, some other types of continuities, namely, uniform continuity, statistical continuity, lacunary statistical continuity as well as theorems related to sequences of functions, and closedness of the set of Abel continuous functions. We also investigate a type of compactness and prove some new results related to this kind of compactness.

\maketitle

\section{Preliminaries}

Throughout this paper, $\Re$  will denote the set of real numbers. We will use boldface letters $\bf{p}$, $\bf{r}$, $\bf{w}$, ... for sequences $\textbf{p}=(p_{n})$, $\textbf{r}=(r_{n})$, $\textbf{w}=(w_{n})$, ... of points in $\Re$. $c$ will denote the set of convergent sequences of points in $\Re$.
A sequence \textbf{p}=$(p_{k})$ of points in $\Re$ is called statistically convergent \cite{Fast} (see also \cite{Fridy}, and \cite{CakalliAstudy}) to an element $\ell$ of $\Re$ if \[\lim_{n \to \infty}\frac{1}{n}|\{k\leq n: |p_{k}-\ell|\geq \varepsilon\}|=0,\] for every $\varepsilon>0,$ and this is denoted by $st-\lim p_{n}=\ell.$\\
A sequence $(p_{k})$ of points in $\Re$  is called lacunarily statistically convergent (\cite{FridyandOrhan1}) to an element $\ell$ of $\Re$  if
\[
\lim_{r\rightarrow\infty}\frac{1}{h_{r}}|\{k\in I_{r}: |p_{k}-\ell| \geq{\varepsilon}\}|=0,
\]
for every $\varepsilon >0$ where $I_{r}=(k_{r-1},k_{r}]$, and $k_{0}=0$, $h_{r}:k_{r}-k_{r-1}\rightarrow \infty$ as $r\rightarrow\infty$ and $\theta=(k_{r})$ is an increasing sequence of positive integers, and this is denoted by $S_{\theta}-\lim_{n\rightarrow\infty}p_{n}=\ell$ (see also \cite{Cakallilacunarysta}). Throughout this paper we assume that $\liminf_{r}\frac{k_{r}}{k_{r-1}}>1$.

A sequence $\textbf{p}=(p_{n})$ of points in $\Re$ is slowly oscillating \cite{Dik} (See also \cite{CakCanDik}, \cite{Cakallislowly}, and \cite{CanakandDik}), denoted by $\textbf{p} \in \textbf{SO}$, if \[\lim_{\lambda \to 1^{+}}\overline{\lim_{n}}\max_{n+1\leq k\leq [\lambda n]}|p_{k}-p_{n}|=0,\] where $[\lambda n]$ denotes the integer part of $\lambda n.$ This is equivalent to the following: $p_{m}-p_{n} \to 0$ whenever $1\leq \frac{m}{n} \to 1$ as $m,n \to \infty.$ In terms of $\varepsilon$ and $\delta,$ this is also equivalent to the case when for any given $\varepsilon>0,$ there exist $\delta=\delta(\varepsilon)>0$ and a positive integer  $N=N(\varepsilon)$ such that $|p_{m}-p_{n}|<\varepsilon$ if $n\geq N(\varepsilon)$ and $n\leq m\leq (1+\delta)n.$

By a method of sequential convergence, or briefly a method, we mean a linear function $G$ defined on a linear subspace of $s$, denoted by $c_{G}$, into $\Re$. A sequence \textbf{p}=$(p_{n})$ is said to be $G$-convergent \cite{ConnorGross} to $\ell$ if \textbf{p}$\in c_{G}$, and $G(\textbf{p})=\ell$. In particular, $\lim$ denotes the limit function $\lim \textbf{p}=\lim_{n}p_{n}$ on the linear space $c$. A method $G$ is called regular if every convergent sequence
\textbf{p}=$(p_{n})$ is $G$-convergent with G(\textbf{p})=$\lim$ \textbf{p}. A method $G$ is called subsequential if whenever \textbf{p} is $G$-convergent with G(\textbf{p})=$\ell,$ then there is a subsequence ($p_{n_{k}}$) of \textbf{p} with $\lim_{k}p_{n_{k}}=\ell.$ A function $f$ is called $G$-continuous (see also \cite{Cakseqdefcomp}, and \cite{Cak7})  if G(f(\textbf{p}))=f(G(\textbf{p})) for any $G$-convergent sequence \textbf{p}. Here we note that for special $G = st -\lim$, $f$ is called statistically continuous \cite{Cakalli4}. For real and complex number sequences, we note that the most important transformation class is
the class of matrix methods. For more information for classical, and modern summability methods see \cite{Boos}.

A sequence $\textbf{p}=(p_{n})$ of real numbers is called Abel
convergent (or Abel summable)  to $\ell$ if
the series $\Sigma_{k=0}^{\infty}p_{k}x^{k}$ is convergent for $0\leq x<1$ and
\[\lim_{x \to 1^{-}}(1-x) \sum_{k=0}^{\infty}p_{k}x^{k}=\ell.\]
In this case we write  $Abel-\lim p_{n}=\ell$. The set of Abel convergent sequences will be denoted by $\textbf{A}$. Abel proved that if $\lim_{n \to \infty}p_{n}=\ell$, then $Abel-\lim p_{n}=\ell$, i.e. every convergent sequence is Abel convergent to the same limit (\cite{Abel}, see also \cite{CVStanojevicandVBStanojevic}, and \cite{FridyandKhan}). As it is known that the converse is not always true in general, as we see that the sequence $((-1)^{n})$ is Abel convergent to $0$ but convergent in the ordinary sense.

\maketitle

\section{ABEL CONTINUITY}

We now give the concept of Abel continuity.
\begin{Def}
A function $f$ is called Abel continuous, denoted by \textbf{f$\in$ AC,} if
it transforms Abel convergent sequences to Abel convergent
sequences, i.e. $(f(p_{n}))$ is Abel convergent to  $f(\ell)$
whenever $(p_{n})$ is Abel convergent to $\ell$.
\end{Def}
We note that this definition of continuity can not be obtained by any $A$-continuity, i.e. by any summability matrix $A$. We also note that sum of two Abel continuous functions is Abel continuous, and composite of two Abel continuous functions is Abel continuous but the product of two Abel continuous functions need not be Abel continuous as it can be seen by considering product of the Abel continuous function $f(t)=t$ with itself.

We note that $G$ defined by $G(\textbf{p})=Abel-lim \textbf{p}$ is a sequential method in the manner of  \cite{Cakseqdefcomp}, but subsequential, so the theorems involving subsequentiality in  \cite{Cakseqdefcomp} cannot be applied to Abel sequential method.

In connection with Abel convergent sequences and convergent sequences the problem arises to investigate the following types of  "continuity" of functions on $\Re$.

\begin{description}
\item[($A$)] $(p_{n}) \in {\textbf{A}} \Rightarrow (f(p_{n})) \in {\textbf{A}}$
\item[$(Ac)$] $(p_{n}) \in {\textbf{A}} \Rightarrow (f(p_{n})) \in {c}$
\item[$(c)$] $(p_{n}) \in {c} \Rightarrow (f(p_{n})) \in {c}$
\item[$(cA)$] $(p_{n}) \in {c} \Rightarrow (f(p_{n})) \in {\textbf{A}}$
\end{description}

We see that $A$ is Abel continuity of $f$, and $(c)$ states the ordinary continuity of $f$. We easily see that (c) implies (cA),  (A) implies (cA), and (Ac) implies (A).  The converses are not always true as the identity function could be taken as a counter example for all the cases.

We note that (c) can be replaced by either statistical continuity, i.e.  $st-\lim f(p_{n})=f(\ell)$  whenever $\textbf{p}=(p_{n})$ is a statistically convergent sequence with $st-\lim p_{n}=\ell$, or lacunary statistical continuity, i.e. $S_{\theta}-\lim_{n\rightarrow \infty} f(p_{n})=f(\ell)$  whenever $\textbf{p}=(p_{n})$ is a lacunary statistically convergent sequence with $S_{\theta}-\lim_{n\rightarrow \infty} p_{n}=\ell$. More generally $(c)$ can be replaced by $G$-sequential continuity of $f$ for any regular subsequential method $G$.

Now we give the implication $(A)$ implies $(c)$, i.e. any Abel continuous function is continuous in the ordinary sense.

\noindent{Theorem 1}
If a function $f$ is Abel continuous on a subset $E$ of $\Re$, then it is continuous on $E$ in the ordinary sense.
\begin{pf}
Suppose that a function $f$ is not continuous on $E$. Then there exists a sequence $(p_{n})$
with $\lim_{n\to \infty } p_{n}=\ell$ such that $(f(p_{n}))$ is not
convergent to $f(\ell)$.
If $(f(p_{n}))$ exists, and $\lim f(p_{n})$ is different from $f(\ell)$ then we
 easily see a contradiction.  Now suppose that $(f(p_{n}))$ has two subsequences of $f(p_{n})$ such that $\lim_{m \to \infty}f(p_{k_{m}})=L_{1}$ and  $\lim_{k \to \infty}f(p_{n_{k}})=L_{2}.$ Since $(p_{n_{k}})$ is subsequence
of $(p_{n}),$  by hypothesis, $\lim_{k \to \infty}f(p_{n_{k}})=f(\ell)$ and $(p_{k_{m}})$ is a subsequence of $(p_{n}),$ by hypothesis $\lim_{m \to \infty}f(p_{k_{m}})=f(\ell).$ This is a contradiction.
If $(f(p_{n}))$ is unbounded above. Then we can find an $n_1$ such that $f(p_{n_{1}})>2^{1}$. There exists a positive integer an $n_{2}>n_{1}$ such that $f(p_{n_{2}})>2^{2}$. Suppose that we have chosen an $n_{k-1}>n_{k-2}$ such that $f(p_{n_{k-1}})>2^{k-1}$. Then we can choose an $n_{k}>n_{k-1}$ such that $f(p_{n_{k}})>2^{k}$. Inductively we can construct a subsequence $(f(p_{n_{k}}))$ of $(f(p_{n}))$ such that $f(p_{n_{k}})>2^{k}$. Since the sequence $(p_{n_{k}})$ is a subsequence of  $(p_{n})$, the subsequence $(p_{n_{k}})$ is convergent so is Abel convergent. But $(f(p_{n_{k}}))$ is not Abel convergent as we see line below. For each positive integer $k$ we have $f(p_{n_{k}})x^k>2^{k}x^k$. The series $\sum_{k=0}^{\infty}2^{k}x^k$ is divergent for any $x$ satisfying $\frac{1}{2}<x<1$, so is the series $\sum_{k=0}^{\infty}f(p_{n_{k}}).x^{k}$. This is a contradiction to the Abel convergence of the sequence $(f(p_{n_{k}}))$. If $(f(p_{n}))$ is unbounded below, similarly $\sum_{k=0}^{\infty}f(p_{k}).x^{k}$ is found to be divergent. The contradiction for  all possible cases to the Abel continuity of $f$ completes the proof of the theorem.
\end{pf}

The converse is not always true for the bounded function $f(t)=\frac{1}{1+t^{2}}$ defined on $\Re$ is an
example. The function $t^{2}$ is another example which is unbounded on $\Re$ as well.

On the other hand not all uniformly continuous functions are Abel continuous. For example the function defined by $f(t)=t^{3}$
is uniformly continuous, but Abel continuous.

\noindent{Corollary 2}
If f is Abel continuous, then it is statistically  continuous.
\begin{pf}
The proof follows from Theorem 1 above, and Corollary 4 in \cite{CakalliAstudy}.
\end{pf}

\noindent{Corollary 3}
If f is Abel continuous, then it is lacunarily statistically sequentially continuous.

\begin{pf}
The proof follows from Theorem 1 above, and Corollary 7 in \cite{Cakallilacunarysta}.
\end{pf}
Now we have the following result.

\noindent{Corollary 4}
If $(p_{n})$ is slowly oscillating, Abel convergent, and $f$ is an Abel continuous function, then $(f(p_{n}))$ is a convergent sequence.

\noindent{Corollary 5}
For any regular subsequential method $G$, any Abel continuous function is $G$-continuous.

For bounded functions we have the following result.

\noindent{Theorem 6} Any bounded Abel continuous function is Cesaro continuous.
\begin{pf}
Let $f$ be a bounded Abel continuous function. Now we are going to obtain that $f$ is Cesaro continuous. To do this take any Cesaro convergent sequence $(p_{n})$ with Cesaro limit $\ell$. Since any Cesaro convergent sequence is Abel convergent to the same value \cite{Rajagopal} (see also \cite{BadiozzamanandThorpeSomebest})
$(p_{n})$ is also Abel convergent to $\ell$. By the assumption that $f$ is Abel continuous, $(f(p_{n}))$ is Abel convergent to $f(\ell)$. By the boundedness of $f$, $(f(p_{n}))$ is bounded. By Corollary to Karamata's Hauptsatz on page 108 in \cite{CVStanojevicandVBStanojevic}, $(f(p_{n}))$ is Cesaro convergent to $f(\ell)$. Thus $f$ is Cesaro continuous at the point $\ell$.  Hence $f$ is continuous at any point in the domain.

\end{pf}

\noindent{Corollary 7} Any bounded Abel continuous function is linear.

\begin{pf}
The proof follows from the preceding theorem, and the theorem on page 73 in \cite{Posner}.
\end{pf}

It is well known that uniform limit of a sequence of continuous functions is continuous. This is also true for Abel continuity, i.e. uniform limit of a sequence of Abel continuous functions is Abel continuous.

\noindent{Theorem 8}
If $(f_{n})$ is a sequence of Abel continuous functions defined on a
subset E of $\Re$ and $(f_{n})$ is uniformly convergent to a
function $f$, then $f$ is Abel continuous on $E$.
\begin{pf} Let $(p_{n})$  be an Abel convergent sequence of real numbers in $E$. Write $Abel-\lim p_{n}=\ell.$ Take any $\varepsilon>0.$  Since $(f_{n})$ is uniformly convergent to  f, there exists a positive integer N such
that
\begin{eqnarray}\label{es1}
  |f_{n}(t)-f(t)|<\varepsilon/3
\end{eqnarray} for all $t\in E$ whenever $n\geq N$. Hence
\begin{eqnarray}\label{es2}
  |(1-x)\sum_{k=0}^{\infty}(f(p_{k})-f_{N}(p_{k}))x^{k}|<\frac{\varepsilon}{3}.
\end{eqnarray}
As $f_{N}$ is Abel continuous, then there exist a $\delta>0$ for
$1-\delta<x<1 $ such that
\begin{eqnarray}\label{es3}
\mid (1-x)
\sum_{k=0}^{\infty}(f_{N}(p_{k})-f_{N}(\ell))x^{k}\mid
<\frac{\varepsilon}{3}.\end{eqnarray}
Now for $1-\delta<x<1,$ it follows
from (\ref{es1}), (\ref{es2}), and (\ref{es3}) that \[\mid
(1-x) \sum_{k=0}^{\infty}(f(p_{k})-f(\ell))x^{k}\mid\leq\mid (1-x)
\sum_{k=0}^{\infty}(f(p_{k})-f_{N}(p_{k}))x^{k}\mid\]\[+\mid (1-x)
\sum_{k=0}^{\infty}(f_{N}(p_{k}))-f_{N}(\ell))x^{k}\mid+\mid (1-x)
\sum_{k=0}^{\infty}(f_{N}(\ell)-f(\ell))x^{k}|\]
\[<
 \frac{\varepsilon}{3}+\frac{\varepsilon}{3}+\frac{\varepsilon}{3}=\varepsilon. \] This
completes the proof of the theorem.
\end{pf}

In the following theorem we prove that the set of Abel continuous functions is a closed subset of the space of continuous functions.

\noindent{Theorem 9}
The set of Abel continuous functions on a subset $E$ of $\Re$ is a closed subset of the set of all continuous functions on $E$, i.e.
$\overline{\textbf{AC}(E)}=\textbf{AC}(E)$, where $\textbf{AC}(E)$ is the set of all Abel continuous functions on $E$, $\overline{\textbf{AC}(E)}$ denotes the set of all cluster points of $\textbf{AC}(E)$.
\begin{pf}
Let f be any element in the closure of $\textbf{AC}(E)$. Then there exists a sequence $(f_{n})$ of points in \textbf{AC}(E) such that $\lim_{n \to \infty}f_{n}=f.$  To show that f is Abel continuous, take any Abel convergent sequence  $(p_{n})$ of points $E$ with Abel limit $\ell$. Let $\varepsilon>0$.
Since $(f_{n})$ is convergent to $f$, there exists a positive integer $N$ such that $$|f_{n}(t)-f(t)|<\varepsilon/6$$ for all $t\in E$ whenever $n\geq N$.
Write \[M=\max \{|f(\ell)-f_{N}(\ell)|, |f(p_{1})-f_{N}(p_{1})|,..., |f(p_{N-1})-f_{N}(p_{N-1})|\},\] and $\delta_{1} =\frac{\varepsilon}{6(N+1)(M+1)}$. Then
we obtain that for any $x$ satisfying $1-\delta_{1}<x<1$ \\
$|(1-x)\sum_{k=0}^{\infty}(f(p_{k})-f_{N}(p_{k}))x^{k}|$\\$\leq|(1-x)\sum_{k=0}^{N-1}(f(p_{k})-f_{N}(p_{k}))x^{k}|+|(1-x)\sum_{k=N}^{\infty}(f(p_{k})-f_{N}(p_{k}))x^{k}|$
\\$<(1-x)MN+|(1-x)\sum_{k=N}^{\infty}(f(p_{k})-f_{N}(p_{k}))x^{k}|<\frac{\varepsilon}{6}+\frac{\varepsilon}{6}=\frac{\varepsilon}{3}$.\\
As $f_{N}$ is Abel continuous, then there exists a $\delta_{2} >0$ such that for $1-\delta_{2}<x<1 $
\[\mid (1-x)
\sum_{k=0}^{\infty}(f_{N}(p_{k})-f_{N}(\ell))x^{k}\mid
<\frac{\varepsilon}{3}.\]
Let $\delta=\min\{\delta_{1},\delta_{2}\}$. Now for $1-\delta<x<1$, we have\\

$|(1-x)\sum_{k=0}^{\infty}f(p_{k})x^{k}-f(\ell)|$
$$\leq\mid (1-x)\sum_{k=0}^{\infty}(f(p_{k})-f_{N}(p_{k}))x^{k}\mid +\mid (1-x) \sum_{k=0}^{\infty}(f_{N}(p_{k}))-f_{N}(\ell))x^{k}\mid\]\[+\mid (1-x)
\sum_{k=0}^{\infty}(f_{N}(\ell)-f(\ell))x^{k}|\]
\[<
 \frac{\varepsilon}{3}+\frac{\varepsilon}{3}+\frac{\varepsilon}{3}=\varepsilon. \]
This completes the proof of the theorem.
\end{pf}

\noindent{Corollary 10}
The set of all Abel continuous functions on a subset $E$ of $\Re$ is a complete subspace of the space of all continuous functions on $E$.

\begin{pf}
The proof follows from the preceding theorem, and the fact that the set of all continuous functions on $E$ is complete.
\end{pf}
Now we can give definition of Abel compactness of a subset of $\Re$.
\begin{Def}
A subset $F$ of $\Re$ is called Abel sequentially compact if whenever
$\textbf{p}=(p_{n})$ is a sequence of point in $F$, there is an Abel convergent
subsequence $\textbf{r}=(r_{k})=(r_{n_{k}})$ of $\textbf{p}$ with
$\lim_{x \to 1^{-}}(1-x)\sum_{k=0}^{\infty}r_{k}x^{k}\in F$.
\end{Def}

\begin{Def}
A real number $\ell$ is said to be in the Abel sequentially closure of a subset
$F$ of $\Re$, denoted by $\overline{F}^{Abel}$ if there is a
sequence $\textbf{p}=(p_{n})$ of points in $F$ such that $Abel-lim
p_{n}=\ell$, and it is called Abel sequentially closed if
$\overline{F}^{Abel}=F$.

It is clear that
$\overline{\phi}^{Abel}=\phi$ and $\overline{\Re}^{Abel}=\Re$.
It is easily seen that $F\subset \overline{F}\subset
\overline{F}^{Abel}$. it is not always true that $\overline
{(\overline{F}^{Abel})}^{Abel}=\overline{F}^{Abel}$; for example,
$\overline {(\overline{\{-1,1\}}^{Abel})}^{Abel}\neq
\overline{\{-1,1\}}^{Abel}$.
\end{Def}

We note that any Abel sequentially closed subset
of Abel sequentially compact subset of $\Re$ is also Abel sequentially compact. Thus
intersection of two Abel sequentially compact, Abel sequentially closed subsets of $\Re$ is
Abel sequentially compact. In general, any intersection of Abel sequentially compact, Abel sequentially closed
subsets of $\Re$ is Abel sequentially compact. A subset of an Abel sequentially compact subset need not to be Abel sequentially compact. For example the interval $]-1,1]$, i.e. the set of real numbers strictly greater than $-1$, and less than or equal to $1$, is a subset of Abel sequentially compact subset $[-1,1]$, i.e. the set of real numbers greater than or equal to $-1$, and less than or equal to $1$, but sequentially Abel compact.
Notice that union of two Abel sequentially compact
subsets of $\Re$ is Abel sequentially compact. We see that any finite union of
Abel sequentially compact subsets of $\Re$ is Abel sequentially compact, but any union of
Abel sequentially compact subsets of $\Re$ is not always Abel sequentially compact.

\noindent{Theorem 11}
Abel continuous image of any Abel sequentially compact subset of $\Re$ is Abel sequentially compact.

\begin{pf}
Although the proof follows from Theorem 7 in \cite{Cakseqdefcomp}, we give a short proof for completeness.
Let $f$ be any Abel continuous function defined on a subset $E$ of $\Re$ and $F$ be any Abel sequentially compact subset of $E$. Take any sequence $\textbf{w}=(w_{n})$ of point in $f(F)$. Write $w_{n}=f(p_{n})$ for each positive integer $n$. Since $F$ is Abel sequentially compact, there exists a Abel convergent subsequence $\textbf{r}=(r_{k})$ of the sequence $\textbf{p}$. Write $Abel-lim\textbf{r}=\ell$. Since $f$ is Abel continuous $Abel-limf(\textbf{r})=f(\ell)$. Thus $f(\textbf{r})=(f(r_{k}))$ is Abel convergent to $f(\ell)$ and a subsequence of the sequence $\textbf{w}$. This completes the proof.
\end{pf}

\noindent{Theorem 12}
If a function $f$ is Abel continuous on a subset $E$ of $\Re$, then $$f(\overline{A}^{Abel})\subset{\overline{(f(A))}^{Abel}}$$ for every subset $A$ of $E$.
\begin{pf}
The proof follows from the regularity of Abel method, and Theorem 8 on page 316 of \cite{Cak7}.
\end{pf}

\noindent{Theorem 13}
For any regular subsequential method $G$, if a subset $F$ of
$\Re$ is $G$-sequentially compact, then it is Abel sequentially compact.
\begin{pf} The proof can be obtained by noticing the regularity
and subsequentiality of $G$ (\cite{Cakseqdefcomp} see $\c{C}akall\i$
for the detail of $G$-sequential compactness).
\end{pf}

\noindent{Corollary 14}
Any compact subset of $\Re$ is Abel sequentially compact.
\begin{pf}
Writing $G=lim\textbf{p}$, and applying the preceding theorem we get the proof.

\end{pf}

\noindent{Corollary 15}
Any Abel sequentially closed subset of an Abel compact sequentially subset of $\Re$ is Abel sequentially compact.
\begin{pf}
Writing $G(\textbf{p})=Abel-lim\textbf{p}$, and applying Theorem 1 on page 597 in \cite{Cakseqdefcomp} we get the proof.
\end{pf}

\noindent{Theorem 16} Any Abel sequentially compact subset of $\Re$ is bounded.
\begin{pf}
Suppose that $E$ is unbounded so that we can construct a sequence $\textbf{p}=(p_{n})$ of points in $E$ such that $p_{n}>2^{n}$, and $p_{n}>p_{n-1}$ for each positive integer $n$. It is easily seen that the sequence $\textbf{p}$ has no Abel convergent subsequence. If it is unbounded below, then similarly we construct a sequence of points in $E$ which has no Abel convergent subsequence. Hence $E$ is not Abel sequentially compact. This contradiction completes the proof of the theorem.
\end{pf}

We note that the converse of this theorem is not always true as it can be seen by considering any open subset of $\Re$.

\noindent{Corollary 17} Any Abel sequentially compact subset of $\Re$ is slowly oscillating compact.

\begin{pf}
The proof follows from Theorem 16.
\end{pf}

\noindent{Corollary 18} Any Abel sequentially compact subset of $\Re$ is $\Delta$-slowly oscillating compact.

\begin{pf}
The proof follows from Theorem 16, and Corollary 17.
\end{pf}

\noindent{Corollary 19} Any Abel sequentially compact subset of $\Re$ is ward compact.

\begin{pf}
The proof follows from Theorem 16 (see \cite{Cakalliforwardcontinuity},and \cite{Cakallionboundedness} for the definition of ward compactness)
\end{pf}

\noindent{Corollary 20} Any Abel sequentially compact subset of $\Re$ is $\delta$-ward compact.

\begin{pf}
The proof follows from Theorem 16 (see \cite{CakallideltaquasiCauchy} for the definition of $\delta$-ward compactness.
\end{pf}

\maketitle

\section{Conclusion}

In this paper we introduce a concept of Abel  continuity, and a concept of Abel sequential compactness, and present theorems related to  this kind of sequential continuity, this kind of sequential compactness, continuity, statistical continuity, lacunary statistical continuity, and uniform continuity. One may expect this investigation to be a useful tool in the field of analysis in modeling various problems occurring in many areas of science, dynamical systems, computer science, information theory, and biological science. On the other hand, we suggest to introduce a concept of fuzzy Abel sequential compactness, and investigate fuzzy Abel continuity for fuzzy functions (see \cite{CakalliPratul} for the definitions and  related concepts in fuzzy setting). However due to the change in settings, the definitions and methods of proofs will not always be the same. We also suggest to investigate a theory in dynamical systems by introducing the following concept: two dynamical systems are called Abel-conjugate if there is one to one, onto, $h$, and $h^{-1}$ are Abel continuous, and $h$ commutes the mappings at each point.

\end{document}